\documentclass[11pt]{amsart}
\usepackage{amsmath,amssymb}
\newtheorem{theorem}{Theorem}
\newtheorem{proposition}[theorem]{Proposition}
\newtheorem{corollary}[theorem]{Corollary}
\newtheorem{lemma}[theorem]{Lemma}

\begin{document}

\title[Asymptotically hyperbolic manifolds]{A volume comparison theorem for asymptotically hyperbolic manifolds}
\author{Simon Brendle and Otis Chodosh}
\address{Department of Mathematics \\ Stanford University \\ 450 Serra Mall, Bldg 380 \\ Stanford, CA 94305} 
\begin{abstract}
We define a notion of renormalized volume of an asymptotically hyperbolic manifold. Moreover, we prove a sharp volume comparison theorem for metrics with scalar curvature at least $-6$. Finally, we show that the inequality is strict unless the metric is isometric to one of the Anti-deSitter-Schwarzschild metrics.
\end{abstract}
\thanks{The first author was supported in part by the National Science Foundation under grant DMS-1201924. The second author was supported in part by a National Science Foundation Graduate Research Fellowship DGE-1147470.}
\maketitle

\section{Introduction}

Let $(\bar{M},\bar{g})$ denote the standard three-dimensional hyperbolic space, so that 
\[\bar{g} = \frac{1}{1 + s^2} \, ds \otimes ds + s^2 \, g_{S^2}.\] 
Let us consider a Riemannian metric $g$ on $M = \bar{M} \setminus K$, where $K$ is a bounded domain with smooth connected boundary. We assume that $g$ is asymptotically hyperbolic in the sense that $|g-\overline{g}|_{\bar{g}} = O(s^{-2-4\delta})$ for some $\delta \in (0,\frac{1}{4})$ and $|\bar{D}(g-\bar{g})|_{\bar{g}} = o(1)$. We define the renormalized volume of $(M,g)$ by 
\[V(M,g) := \lim_{i \to \infty} (\text{\rm vol}(\Omega_i \cap M,g) - \text{\rm vol}(\Omega_i,\bar{g})),\] 
where $\Omega_i$ is an arbitrary exhaustion of $\bar{M}$ by compact sets. The condition $|g-\overline{g}|_{\bar{g}} = O(s^{-2-4\delta})$ guarantees that the quantity $V(M,g)$ does not depend on the choice of the exhaustion $\Omega_i$. Clearly, $V(\bar{M},\bar{g}) = 0$.

As an example, let us consider the Anti-deSitter-Schwarzschild manifold with mass $m > 0$. To that end, let $s_0 = s_0(m)$ denote the unique positive solution of the equation $1 + s_0^2 - m \, s_0^{-1} = 0$. We then consider the manifold $\bar{M}_m = \bar{M} \setminus \{s \leq s_0(m)\}$ equipped with the Riemannian metric
\[\bar{g}_m = \frac{1}{1 + s^2 - m \, s^{-1}} \, ds \otimes ds + s^2 \, g_{S^2}.\]
The boundary $S^2 \times \{s_0(m)\}$ is an outermost minimal surface, which is referred to as the horizon. Moreover, it is easy to see that $|\bar{g}_m-\bar{g}|_{\bar{g}} = O(s^{-3})$, so $\bar{g}$ satisfies the asymptotic assumptions above. The renormalized volume of $(\bar{M}_m,\bar{g}_m)$ is given by 
\[V(\bar{M}_m,\bar{g}_m) = \lim_{r \to \infty} \bigg ( \int_{s_0(m)}^r \frac{4\pi s^2}{\sqrt{1+s^2-m \, s^{-1}}} \, ds - \int_0^r \frac{4\pi s^2}{\sqrt{1+s^2}} \, ds \bigg ).\]
We now state the main result of this paper.

\begin{theorem} 
\label{main.theorem}
Let us consider a Riemannian metric $g$ on $M = \bar{M} \setminus K$, where $K$ is a compact set with smooth connected boundary. We assume that $g$ has the following properties: 
\begin{itemize} 
\item The manifold $(M,g)$ is asymptotically hyperbolic in the sense that $|g-\overline{g}|_{\bar{g}} = O(s^{-2-4\delta})$ and $|\bar{D}(g-\bar{g})|_{\bar{g}} = o(1)$.
\item The scalar curvature of $g$ is at least $-6$.
\item The boundary $\partial M$ is an outermost minimal surface with respect to $g$, and we have $\text{\rm area}(\partial M,g) \geq \text{\rm area}(\partial \bar{M}_m,\bar{g}_m)$ for some $m>0$.
\end{itemize}
Then $V(M,g) \geq V(\bar{M},\bar{g}_m)$. Moreover, if equality holds, then $g$ is isometric to $\bar{g}_m$.
\end{theorem}

We note that our asymptotic assumptions are quite weak: in particular, $g$ and $\bar{g}_m$ may have different mass at infinity. An immediate consequence of Theorem \ref{main.theorem} is that the function $m \mapsto V(\bar{M}_m,\bar{g}_m)$ is strictly monotone increasing. This fact is not obvious, as $s_0(m)$ is an increasing function of $m$.

Theorem \ref{main.theorem} is motivated in part by Bray's volume comparison theorem \cite{Bray1} for three-manifolds with scalar curvature at least $6$, as well as by a rigidity result due to Llarull \cite{Llarull}. A survey of this and other rigidity results involving scalar curvature can be found in \cite{Brendle1}.

The proof of Theorem \ref{main.theorem} uses two main ingredients. The first is the weak inverse mean curvature flow, which was introduced in the ground-breaking work of Huisken and Ilmanen \cite{Huisken-Ilmanen} on the Riemannian Penrose inequality (see also \cite{Bray2}, where an alternative proof is given). The inverse mean curvature flow has also been considered as a possible tool for proving a version of the Penrose inequality for asymptotically hyperbolic manifolds; see \cite{Neves}, \cite{Wang}. More recently, the inverse mean curvature flow was used in \cite{Brendle-Hung-Wang} to prove a sharp Minkowski-type inequality for surfaces in the Anti-deSitter-Schwarzschild manifold.

The second ingredient in our argument is an isoperimetric principle which asserts that a coordinate sphere in the standard Anti-deSitter-Schwarzschild manifold has smallest area among all surfaces that are homologous to the horizon and enclose the same amount of volume. This inequality was established in \cite{Corvino-Gerek-Greenberg-Krummel} following earlier work by Bray \cite{Bray1}. In fact, it is known that the coordinate spheres are the only embedded hypersurfaces with constant mean curvature in the Anti-deSitter-Schwarzschild manifold (see \cite{Brendle2}).

Our approach also shares common features with a result of Bray and Miao \cite{Bray-Miao}, which gives a sharp bound for the capacity of a surface in a three-manifold of nonnegative scalar curvature. 

\section{Proof of Theorem \ref{main.theorem}}

Let $(M,g)$ be a Riemannian manifold which satisfies the assumptions of Theorem \ref{main.theorem}, and let $(\bar{M}_m,\bar{g}_m)$ be an Anti-deSitter-Schwarzschild manifold satisfying $\text{\rm area}(\partial M,g) \geq \text{\rm area}(\partial \bar{M}_m,\bar{g}_m)$. For abbreviation, let $A = \text{\rm area}(\partial M,g)$ and $\bar{A} = \text{\rm area}(\partial \bar{M}_m,\bar{g}_m)$.

Let $\Sigma_t$ denote the weak solution of the inverse mean curvature flow in $(M,g)$ with the initial surface $\Sigma_0 = \partial M$. For each $t$, we denote by $\Omega_t \subset \bar{M}$ the region bounded by $\Sigma_t$. 

\begin{proposition}
\label{barrier}
Let $\delta \in (0,\frac{1}{4})$ be as above. Then we have $\{s \leq e^{\frac{(1-\delta)t}{2}}\} \subset \Omega_t$ for $t$ sufficiently large.
\end{proposition}

\textbf{Proof.} 
The coordinate sphere $S^2 \times \{s\}$ has mean curvature $2 + o(1)$ for $s$ large. Hence, we can find a real number $t_0$ such that the surfaces $S_t = \{s=e^{\frac{(1-2\delta)t}{2}}\}$ move with a speed less than $\frac{1}{H}$ for $t \geq t_0$. By the Weak Existence Theorem 3.1 in \cite{Huisken-Ilmanen}, the regions $\Omega_t$ will eventually contain every given compact set. Hence, we can find a real number $\tau$ such that $\{s \leq e^{\frac{(1-2\delta)t_0}{2}}\} \subset \Omega_\tau$. By the maximum principle (cf. Theorem 2.2 in \cite{Huisken-Ilmanen}), we have $\{s \leq e^{\frac{(1-2\delta)(t-\tau+t_0)}{2}}\} \subset \Omega_t$ for $t \geq \tau$. From this, the assertion follows. \\

Since the boundary $\partial M$ is an outermost minimal surface, we have $\text{\rm area}(\Sigma_t,g) = e^t \, A$. Moreover, it is well known that the quantity
\[m_H(\Sigma_t) = \text{\rm area}(\Sigma_t,g)^{\frac{1}{2}} \, \bigg ( 16\pi - \int_{\Sigma_t} (H_g^2-4) \, d\mu_g \bigg )\] 
is monotone increasing in $t$. 

\begin{proposition} 
\label{volume.estimate}
For each $\tau \geq 0$, we have 
\[\text{\rm vol}(\Omega_\tau \cap M,g) \geq \int_0^\tau e^{\frac{3t}{2}} \, A^{\frac{3}{2}} \, (4 \, e^t \, A + 16\pi - e^{-\frac{t}{2}} \, A^{-\frac{1}{2}} \, m_H(\Sigma_t))^{-\frac{1}{2}} \, dt.\] 
\end{proposition} 

\textbf{Proof.} 
By the co-area formula, we have 
\[\int_{\Omega_\tau \cap M} \psi \, H_g \, d\text{\rm vol}_g = \int_0^\tau \bigg ( \int_{\Sigma_t} \psi \, d\mu_g \bigg ) \, dt\] 
for every nonnegative measurable function $\psi$. Hence, if we put 
\[\psi = \begin{cases} \frac{1}{H_g} & \text{\rm if $H_g > 0$} \\ \infty & \text{\rm if $H_g=0$,} \end{cases}\] 
then we obtain 
\[\text{\rm vol}(\Omega_\tau \cap M,g) \geq \int_0^\tau \bigg ( \int_{\Sigma_t} \psi \, d\mu_g \bigg ) \, dt.\] 
Moreover, it follows from H\"older's inequality that 
\begin{align*} 
\int_{\Sigma_t} \psi \, d\mu_g 
&\geq \text{\rm area}(\Sigma_t,g)^{\frac{3}{2}} \, \bigg ( \int_{\Sigma_t} H_g^2 \, d\mu_g \bigg )^{-\frac{1}{2}} \\ 
&= \text{\rm area}(\Sigma_t,g)^{\frac{3}{2}} \, (4 \, \text{\rm area}(\Sigma_t,g) + 16\pi - \text{\rm area}(\Sigma_t,g)^{-\frac{1}{2}} \, m_H(\Sigma_t))^{-\frac{1}{2}} \\ 
&= e^{\frac{3t}{2}} \, A^{\frac{3}{2}} \, (4 \, e^t \, A + 16\pi - e^{-\frac{t}{2}} \, A^{-\frac{1}{2}} \, m_H(\Sigma_t))^{-\frac{1}{2}}. 
\end{align*} 
Putting these facts together, the assertion follows. \\

\begin{corollary} 
\label{volume.estimate.2}
We have 
\[2 \, \text{\rm vol}(\Omega_\tau \cap M,g) \geq \int_0^\tau e^t \, A^{\frac{3}{2}} \, ((1-e^{-\frac{3t}{2}}) \, A + 4\pi \, (e^{-t} - e^{-\frac{3t}{2}}))^{-\frac{1}{2}} \, dt.\] 
\end{corollary} 

\textbf{Proof.}
Using the monotonicity of $m_H(\Sigma_t)$, we obtain 
\[m_H(\Sigma_t) \geq m_H(\Sigma_0) = 4 \, A^{\frac{1}{2}} \, (A + 4\pi).\] 
This implies 
\[2 \, \text{\rm vol}(\Omega_\tau \cap M,g) \geq \int_0^\tau e^{\frac{3t}{2}} \, A^{\frac{3}{2}} \, ((e^t-e^{-\frac{t}{2}}) \, A + 4\pi \, (1-e^{-\frac{t}{2}}))^{-\frac{1}{2}} \, dt.\] From this, the assertion follows. \\

\begin{proposition} 
\label{iso}
Let $\Omega$ be a domain in $\bar{M}$ such that $\{s \leq s_0(m)\} \subset \Omega$, and let $\Sigma$ denote the boundary of $\Omega$. Then 
\[2 \, \text{\rm vol}(\Omega \cap \bar{M}_m,\bar{g}_m) \leq \int_0^{\bar{\tau}} e^t \, \bar{A}^{\frac{3}{2}} \, ((1-e^{-\frac{3t}{2}}) \, \bar{A} + 4\pi \, (e^{-t} - e^{-\frac{3t}{2}}))^{-\frac{1}{2}} \, dt,\] 
where $\bar{\tau}$ is defined by $\text{\rm area}(\Sigma,\bar{g}_m) = e^{\bar{\tau}} \, \bar{A}$.
\end{proposition} 

\textbf{Proof.} 
If $\Sigma$ is a coordinate sphere in $(\bar{M}_m,\bar{g}_m)$, then we have 
\[2 \, \text{\rm vol}(\Omega \cap \bar{M}_m,\bar{g}_m) = \int_0^{\bar{\tau}} e^t \, \bar{A}^{\frac{3}{2}} \, ((1-e^{-\frac{3t}{2}}) \, \bar{A} + 4\pi \, (e^{-t} - e^{-\frac{3t}{2}}))^{-\frac{1}{2}} \, dt,\] 
where $\bar{\tau}$ is defined by $\text{\rm area}(\Sigma,\bar{g}_m) = e^{\bar{\tau}} \, \bar{A}$. On the other hand, it is known (cf. \cite{Corvino-Gerek-Greenberg-Krummel}, Theorem 4.2) that the coordinate spheres in $(\bar{M}_m,\bar{g}_m)$ enclose the largest volume for any given surface area. Putting these facts together, the assertion follows. \\

Let us consider a sequence of times $\tau_i \to \infty$. Moreover, we define a sequence of times $\bar{\tau}_i \to \infty$ by $\text{\rm area}(\Sigma_{\tau_i},\bar{g}_m) = e^{\bar{\tau}_i} \, \bar{A}$. By Proposition \ref{barrier}, we have $s \geq e^{\frac{(1-\delta)t}{2}}$ on $\Sigma_t$ if $t$ is large enough. This implies 
\[|g-\bar{g}_m|_{\bar{g}_m} \leq O(s^{-2-4\delta}) \leq O(e^{-(1-\delta)(1+2\delta)t})\] 
at each point on $\Sigma_t$. From this, we deduce that 
\begin{align*} 
e^{\tau_i} \, A 
&= \text{\rm area}(\Sigma_{\tau_i},g) \\ 
&= \text{\rm area}(\Sigma_{\tau_i},\bar{g}_m) \, (1 + O(e^{-(1-\delta)(1+2\delta)\tau_i})) \\ 
&= e^{\bar{\tau}_i} \, \bar{A} \, (1 + O(e^{-(1-\delta)(1+2\delta)\tau_i})). 
\end{align*} 
Thus, we conclude that 
\[\tau_i = \bar{\tau}_i - \alpha + O(e^{-(1-\delta)(1+2\delta)\tau_i}),\] 
where $\alpha = \log (A/\bar{A}) \geq 0$. Note that $(1-\delta)(1+2\delta) > 1$ since $\delta \in (0,\frac{1}{4})$.

By Corollary \ref{volume.estimate.2}, we have 
\begin{align*} 
2 \, \text{\rm vol}(\Omega_{\tau_i} \cap M,g) 
&\geq \int_0^{\tau_i} e^t \, A^{\frac{3}{2}} \, ((1-e^{-\frac{3t}{2}}) \, A + 4\pi \, (e^{-t} - e^{-\frac{3t}{2}}))^{-\frac{1}{2}} \, dt \\ 
&= \int_\alpha^{\tau_i+\alpha} e^{t-\alpha} \, A^{\frac{3}{2}} \, ((1-e^{-\frac{3t-3\alpha}{2}}) \, A + 4\pi \, (e^{-t+\alpha} - e^{-\frac{3t-3\alpha}{2}}))^{-\frac{1}{2}} \, dt \\ 
&= \int_\alpha^{\tau_i+\alpha} e^t \, \bar{A}^{\frac{3}{2}} \, ((1-e^{-\frac{3t-3\alpha}{2}}) \, \bar{A} + 4\pi \, (e^{-t} - e^{-\frac{3t-\alpha}{2}}))^{-\frac{1}{2}} \, dt. 
\end{align*} 
On the other hand, we have 
\[2 \, \text{\rm vol}(\Omega_{\tau_i} \cap \bar{M}_m,\bar{g}_m) \leq \int_0^{\bar{\tau}_i} e^t \, \bar{A}^{\frac{3}{2}} \, ((1-e^{-\frac{3t}{2}}) \, \bar{A} + 4\pi \, (e^{-t} - e^{-\frac{3t}{2}}))^{-\frac{1}{2}} \, dt\] 
by Proposition \ref{iso}. Putting these facts together, we obtain 
\begin{align*} 
&2 \, (V(M,g) - V(\bar{M}_m,\bar{g}_m)) \\ 
&= \limsup_{i \to \infty} 2 \, (\text{\rm vol}(\Omega_{\tau_i} \cap M,g) - \text{\rm vol}(\Omega_{\tau_i} \cap \bar{M}_m,\bar{g}_m)) \\ 
&\geq \limsup_{i \to \infty} \bigg ( \int_\alpha^{\tau_i+\alpha} e^t \, \bar{A}^{\frac{3}{2}} \, ((1-e^{-\frac{3t-3\alpha}{2}}) \, \bar{A} + 4\pi \, (e^{-t} - e^{-\frac{3t-\alpha}{2}}))^{-\frac{1}{2}} \, dt \\ 
&\hspace{20mm} - \int_0^{\bar{\tau}_i} e^t \, \bar{A}^{\frac{3}{2}} \, ((1-e^{-\frac{3t}{2}}) \, \bar{A} + 4\pi \, (e^{-t} - e^{-\frac{3t}{2}}))^{-\frac{1}{2}} \, dt \bigg ) \\ 
&= \limsup_{i \to \infty} \bigg ( \int_\alpha^{\bar{\tau}_i} e^t \, \bar{A}^{\frac{3}{2}} \, ((1-e^{-\frac{3t-3\alpha}{2}}) \, \bar{A} + 4\pi \, (e^{-t} - e^{-\frac{3t-\alpha}{2}}))^{-\frac{1}{2}} \, dt \\ 
&\hspace{20mm} - \int_0^{\bar{\tau}_i} e^t \, \bar{A}^{\frac{3}{2}} \, ((1-e^{-\frac{3t}{2}}) \, \bar{A} + 4\pi \, (e^{-t} - e^{-\frac{3t}{2}}))^{-\frac{1}{2}} \, dt \bigg ) \\ 
&= \bar{A}^{\frac{3}{2}} \, I(\alpha), 
\end{align*} 
where 
\begin{align*} 
I(\alpha) &= \int_\alpha^\infty e^t \, \Big [ ((1-e^{-\frac{3t-3\alpha}{2}}) \, \bar{A} + 4\pi \, (e^{-t} - e^{-\frac{3t-\alpha}{2}}))^{-\frac{1}{2}} \\ 
&\hspace{20mm} - ((1-e^{-\frac{3t}{2}}) \, \bar{A} + 4\pi \, (e^{-t} - e^{-\frac{3t}{2}}))^{-\frac{1}{2}} \Big ] \, dt \\ 
&- \int_0^\alpha e^t \, ((1-e^{-\frac{3t}{2}}) \, \bar{A} + 4\pi \, (e^{-t} - e^{-\frac{3t}{2}}))^{-\frac{1}{2}} \, dt. 
\end{align*} 
It is shown in the appendix that the function $I(\alpha)$ is positive for all $\alpha>0$. Thus, we conclude that $V(M,g) \geq V(\bar{M}_m,\bar{g}_m)$.

Finally, we analyze the case of equality. Suppose that $V(M,g) = V(\bar{M}_m,\bar{g}_m)$. Then $I(\alpha) \leq 0$, which implies that $\alpha=0$. Moreover, the difference 
\[2 \, \text{\rm vol}(\Omega_{\tau_i} \cap M,g) - \int_0^{\tau_i} e^t \, A^{\frac{3}{2}} \, ((1-e^{-\frac{3t}{2}}) \, A + 4\pi \, (e^{-t} - e^{-\frac{3t}{2}}))^{-\frac{1}{2}} \, dt\] 
must converge to $0$ as $i \to \infty$. Using Proposition \ref{volume.estimate}, we conclude that $m_H(\Sigma_t) = m_H(\Sigma_0)$ for all $t$. This implies that $g$ is the isometric to one of the standard Anti-deSitter-Schwarzschild metrics. Since $\alpha=0$, the manifolds $(M,g)$ and $(\bar{M}_m,\bar{g}_m)$ have the same boundary area. Therefore, they are isometric.

\appendix

\section{Positivity of the function $I(\alpha)$}

In this section, we show that $I(\alpha)>0$ for all $\alpha>0$. We begin with a lemma: 

\begin{lemma} 
\label{aux}
Let $\varepsilon$ and $\mu$ be positive real numbers. If the ratio $\frac{\varepsilon}{\mu}$ is sufficiently small, then we have 
\[3\mu \int_0^\infty e^{-\frac{t}{2}} \, (\varepsilon+(1-e^{-\frac{3t}{2}}) \mu)^{-\frac{3}{2}} \, dt \geq 4 \, \varepsilon^{-\frac{1}{2}} + \mu^{-\frac{1}{2}}.\] 
\end{lemma}

\textbf{Proof.} 
It is elementary to check that 
\[e^t \geq 1 + \frac{2}{3} \, (1-e^{-\frac{3t}{2}}),\] 
hence 
\[e^{-\frac{t}{2}} \geq e^{-\frac{3t}{2}} + \frac{2}{3} \, e^{-\frac{3t}{2}} \, (1-e^{-\frac{3t}{2}})\] 
for all $t \geq 0$. This implies 
\begin{align*} 
&\int_0^1 e^{-\frac{t}{2}} \, (\varepsilon+1-e^{-\frac{3t}{2}})^{-\frac{3}{2}} \, dt \\ 
&\geq \int_0^1 e^{-\frac{3t}{2}} \, (\varepsilon+1-e^{-\frac{3t}{2}})^{-\frac{3}{2}} \, dt + \frac{2}{3} \int_0^1 e^{-\frac{3t}{2}} \, (1-e^{-\frac{3t}{2}}) \, (\varepsilon+1-e^{-\frac{3t}{2}})^{-\frac{3}{2}} \, dt \\ 
&= \int_0^1 e^{-\frac{3t}{2}} \, (\varepsilon+1-e^{-\frac{3t}{2}})^{-\frac{3}{2}} \, dt +
\frac{2}{3} \int_0^1 e^{-\frac{3t}{2}} \, (1-e^{-\frac{3t}{2}})^{-\frac{1}{2}} \, dt - o(1) \\ 
&= \frac{4}{3} \, \varepsilon^{-\frac{1}{2}} - \frac{4}{3} \, (\varepsilon+1-e^{-\frac{3}{2}})^{-\frac{1}{2}} + \frac{8}{9} \, (1-e^{-\frac{3}{2}})^{\frac{1}{2}} - o(1) 
\end{align*} 
for $\varepsilon > 0$ sufficiently small. Hence, we obtain 
\begin{align*} 
&\int_0^\infty e^{-\frac{t}{2}} \, (\varepsilon+1-e^{-\frac{3t}{2}})^{-\frac{3}{2}} \, dt \\ 
&\geq \frac{4}{3} \, \varepsilon^{-\frac{1}{2}} - \frac{4}{3} \, (\varepsilon+1-e^{-\frac{3}{2}})^{-\frac{1}{2}} + \frac{8}{9} \, (1-e^{-\frac{3}{2}})^{\frac{1}{2}} + (\varepsilon+1)^{-\frac{3}{2}} \int_1^\infty e^{-\frac{t}{2}} \, dt - o(1) \\ 
&\geq \frac{4}{3} \, \varepsilon^{-\frac{1}{2}} - \frac{4}{3} \, (1-e^{-\frac{3}{2}})^{-\frac{1}{2}} + \frac{8}{9} \, (1-e^{-\frac{3}{2}})^{\frac{1}{2}} + 2 \, e^{-\frac{1}{2}} - o(1) \\ 
&\geq \frac{4}{3} \, \varepsilon^{-\frac{1}{2}} + \frac{1}{3} 
\end{align*} 
if $\varepsilon>0$ is small enough. This proves the assertion for $\mu=1$. The general case follows by scaling. \\

We now consider the function 
\begin{align*} 
I_\varepsilon(\alpha) &= \int_\alpha^\infty e^t \, \Big [ (\varepsilon+(1-e^{-\frac{3t-3\alpha}{2}}) \, \bar{A} + 4\pi \, (e^{-t} - e^{-\frac{3t-\alpha}{2}}))^{-\frac{1}{2}} \\ 
&\hspace{20mm} - (\varepsilon+(1-e^{-\frac{3t}{2}}) \, \bar{A} + 4\pi \, (e^{-t} - e^{-\frac{3t}{2}}))^{-\frac{1}{2}} \Big ] \, dt \\ 
&- \int_0^\alpha e^t \, (\varepsilon+(1-e^{-\frac{3t}{2}}) \, \bar{A} + 4\pi \, (e^{-t} - e^{-\frac{3t}{2}}))^{-\frac{1}{2}} \, dt. 
\end{align*} 
Then 
\begin{align*} 
\frac{d}{d\alpha} I_\varepsilon(\alpha) &= \int_\alpha^\infty e^t \, \frac{d}{d\alpha} \Big [ (\varepsilon+(1-e^{-\frac{3t-3\alpha}{2}}) \, \bar{A} + 4\pi \, (e^{-t} - e^{-\frac{3t-\alpha}{2}}))^{-\frac{1}{2}} \Big ] \, dt - e^\alpha \, \varepsilon^{-\frac{1}{2}} \\ 
&= \frac{1}{4} \, (3 \, e^{\frac{3\alpha}{2}} \, \bar{A} + 4\pi \, e^{\frac{\alpha}{2}}) \\ 
&\hspace{5mm} \cdot \int_\alpha^\infty e^{-\frac{t}{2}} \, (\varepsilon+(1-e^{-\frac{3t-3\alpha}{2}}) \, \bar{A} + 4\pi \, (e^{-t} - e^{-\frac{3t-\alpha}{2}}))^{-\frac{3}{2}} \, dt - e^\alpha \, \varepsilon^{-\frac{1}{2}} \\ 
&= \frac{e^\alpha}{4} \, (3 \, \bar{A} + 4\pi \, e^{-\alpha}) \\ 
&\hspace{5mm} \cdot \int_0^\infty e^{-\frac{t}{2}} \, (\varepsilon+(1-e^{-\frac{3t}{2}}) \, \bar{A} + 4\pi \, e^{-\alpha} \, (e^{-t} - e^{-\frac{3t}{2}}))^{-\frac{3}{2}} \, dt - e^\alpha \, \varepsilon^{-\frac{1}{2}} \\ 
&\geq \frac{e^\alpha}{4} \, (3 \, \bar{A} + 4\pi \, e^{-\alpha}) \\ 
&\hspace{5mm} \cdot \int_0^\infty e^{-\frac{t}{2}} \, \Big ( \varepsilon+(1-e^{-\frac{3t}{2}}) \, (\bar{A} + \frac{4\pi}{3} \, e^{-\alpha}) \Big )^{-\frac{3}{2}} \, dt - e^\alpha \, \varepsilon^{-\frac{1}{2}}, 
\end{align*} 
where in the last step we have used the inequality $e^{-t} - e^{-\frac{3t}{2}} \leq \frac{1}{3} \, (1-e^{-\frac{3t}{2}})$. Hence, if the ratio $\frac{\varepsilon}{\bar{A} + \frac{4\pi}{3} \, e^{-\alpha}}$ is sufficiently small, then 
\[\frac{d}{d\alpha} I_\varepsilon(\alpha) \geq \frac{e^\alpha}{4} \, (\bar{A} + \frac{4\pi}{3} \, e^{-\alpha})^{-\frac{1}{2}}\] 
by Lemma \ref{aux}. Since $I(\alpha) = \lim_{\varepsilon \to 0} I_\varepsilon(\alpha)$ for each $\alpha \geq 0$, we conclude that the function $I(\alpha)$ is strictly monotone increasing. In particular, $I(\alpha) > 0$ for all $\alpha>0$.

\end{document}